\documentclass{article}
\usepackage[utf8]{inputenc}
\usepackage{authblk}
\usepackage[colorlinks]{hyperref}

\title{An Exact Integral-to-Sum Relation for Products of Bessel Functions}

\author[1,2,*]{Oliver H.E. Philcox}
\author[3,4]{Zachary Slepian}
\affil[1]{\footnotesize Department of Astrophysical Sciences, Princeton University, Princeton, NJ 08540, USA}
\affil[2]{\footnotesize School of Natural Sciences, Institute for Advanced Study, 1 Einstein Drive, Princeton, NJ 08540, USA}
\affil[3]{\footnotesize Department of Astronomy, University of Florida, 211 Bryant Space Science Center, Gainesville, FL 32611, USA}
\affil[4]{\footnotesize Physics Division, Lawrence Berkeley National Laboratory, 1 Cyclotron Road, Berkeley, CA 94709, USA}
\affil[*]{Electronic Address: \href{mailto: ohep2@cantab.ac.uk}{ohep2@cantab.ac.uk}}

\usepackage[a4paper, total={6in, 8in}]{geometry}
\usepackage{amsmath}
\usepackage{amssymb}

\usepackage[numbers]{natbib}
\usepackage{graphicx}
\usepackage{subcaption}
\renewcommand{\vec}[1]{\boldsymbol{#1}}

\newcommand{\vk}{\vec k}

\newcommand{\vr}{\vec r}

\newcommand{\hr}{\hat{\vec r}}

\usepackage{color}
\usepackage{enumitem}
\def\beq{\begin{eqnarray}}
\def\eeq{\end{eqnarray}}

\definecolor{darkgreen}{RGB}{0,120,0}

\newcommand{\new}[1]{#1}
\newcommand{\neww}[1]{#1}

\begin{document}

\maketitle

\begin{abstract}
        A useful identity relating the infinite sum of two Bessel functions to their infinite integral was discovered in Dominici \textit{et al.} (2012). Here, we extend this result to products of $N$ Bessel functions, and show it can be straightforwardly proven using the Abel-Plana theorem, \neww{or the Poisson summation formula}. For $N=2$, the proof is much simpler than that of \neww{Dominici \textit{et al.}}, and significantly enlarges the \new{range of validity}.
\end{abstract}

\section{Introduction}\label{sec: intro}
Integrals of Bessel functions appear in many guises across the physical sciences, populating fields as diverse as atomic physics, \new{classical mechanics} and cosmology \new{\citep[e.g.,][]{1993PhRvA..48.1134B,1995PhRvC..51.2031M,2007PhRvD..76h3523F,2009PhRvD..79f3502H,2010JOSAA..27.2144B,2012PhRvC..86a4003W,2012CQGra..29v4007A,2013MNRAS.434.1808M,2015LRR....18....1C,2015OExpr..2333044F,2015PhRvD..91b3508V,2016JOSAB..33C..30S,2016TMP...188.1197K,2017JCAP...11..039F,2019MNRAS.483.2078Y,2019LRR....22....6H,2018arXiv181202728S,2021MNRAS.501.4004P}.} Of particular importance is their appearance in the basis decomposition of spherically- and circularly-symmetric functions. As a concrete example, consider the inverse Fourier transform of an isotropic function $\tilde f(\vk)\equiv \tilde f(k)$ in three dimensions, defined by
\beq\label{eq: f-intro}
    f(\vr) &=& \int \frac{d^3k}{(2\pi)^3}\,e^{i\vk\cdot\vr}\tilde f(k).
\eeq
Via the well-known plane-wave identity \citep[e.g.,][Eq.\,16.63]{arfken2013mathematical}, the exponential may be expanded as a sum of spherical Bessel functions of the first kind, $j_\ell(kr)$, and the angular part of \new{the integral \eqref{eq: f-intro} can be} performed analytically, leading to
\beq
    f(\vr) &\equiv& f(r) = \int_0^\infty\frac{k^2dk}{2\pi^2}\,j_{\neww{0}}(kr)\tilde f(k).
\eeq
If instead one wishes to evaluate the function $f$ at the \textit{difference} of two positions, $\vr-\vr'$, the same \new{identity} can be used to show
\beq
    f(\vr-\vr') &=& \sum_{\ell=0}^\infty (2\ell+1)\mathcal{L}_\ell(\hr\cdot\hr')\int_0^\infty \frac{k^2dk}{2\pi^2}\,j_\ell(kr)j_{\ell}(kr')\tilde f(k),
\eeq
where \new{$\mathcal{L}_\ell(\hr\cdot\hr')$} is a Legendre polynomial of order $\ell$. This now involves the infinite integral of \textit{two} spherical Bessel functions, \new{which form the coefficients of a} Legendre series in the separation angle $\hr\cdot\hr'$ \citep[e.g.,][]{2019arXiv191200065S}. Whilst not strictly the subject of this work, the above examples serve to illustrate the ubiquity of (spherical) Bessel function integrals.

Given this, much work has been devoted to numerical computation of Bessel function integrals. Of particular note is the FFTLog algorithm \citep{hamilton2000} in which one expands a general function (subject to a set of regularity conditions) as a complex power law, allowing analytic computations of its Bessel-weighted infinite integral. This is both an accurate and fast procedure, and has found great use in the field of cosmology \citep[e.g.,][]{simonovicfftlog}. 

Of interest to this work are the results of \citep{Dominici_2012}, which demonstrated the following relation between the infinite integral of two Bessel functions and their sum:
\beq\label{eq: dominici-simple}
    \int_0^\infty dt\,\frac{J_\nu(at)
    J_\nu(bt)}{t} = \sum_{m=0}^\infty \varepsilon_m\frac{J_\nu(am)J_\nu(bm)}{m}
\eeq
where $J_\nu(x)$ is a Bessel function of the first kind, $\nu$ is a positive half-integer, $a,b\in[0,\pi]$ and
\beq\label{eq: epsilon-m-def}
    \varepsilon_m = \begin{cases} 1/2 & m=0\\ 1 & m\geq 1. \end{cases}
\eeq
The integral \eqref{eq: dominici-simple} may be additionally expressed in terms of spherical Bessel functions of the first kind, giving
\beq
    \int_0^\infty dt\,j_{\ell}(at)j_{\ell}(bt) = \sum_{m=0}^\infty \varepsilon_m\,j_{\ell}(am)j_{\ell}(bm)
\eeq
for non-negative integer \new{$\ell$}. \eqref{eq: dominici-simple} has been used in a variety of contexts \citep{2012JChPh.136j4102L,2016RSPSA.47260421F,2018arXiv181202728S,2019PhRvL.123b0201M,2020JCAP...04..011M}, and was proven in \citep{Dominici_2012} as a special case of the more general relation
\beq
    \int_0^\infty dt\,\frac{J_\mu(at)
    J_\nu(bt)}{t^{\mu+\nu-2k}} = \sum_{m=0}^\infty \varepsilon_m\frac{J_\mu(am)J_\nu(bm)}{m^{\mu+\nu-2k}},
\eeq
for $0<b<a<\pi$, $\mathrm{Re}(\mu)>2k-1/2$, $\mathrm{Re}(\nu)>-1/2$, and $k\in\mathbb{N}_0$ (where $\mathbb{N}_0$ is the set of all natural numbers and zero). To prove this, the authors of \citep{Dominici_2012} use the known result for this infinite integral in terms of Gauss' hypergeometric function (e.g.\,\citep{1994tisp.book.....G}), then manipulate this solution. Below, we show that a simpler proof is possible via the Abel-Plana theorem.

\section{Main Result}\label{sec: main-result}
The principal new result of this work is the following:
\beq\label{eq: main-result}
    \boxed{\int_{0}^\infty dt\, t^{2k}\prod_{j=1}^N\left[t^{-\nu_j}J_{\nu_j}(a_jt)\right] = \sum_{m=0}^\infty \varepsilon_m\, m^{2k}\prod_{j=1}^N\left[m^{-\nu_j}J_{\nu_j}(a_jm)\right]}
\eeq
\textit{i.e.}, that the integral of $N$ \new{Bessel functions} can be written as an infinite sum for arbitrary $N>0$. This uses the definition \eqref{eq: epsilon-m-def} and \new{is valid for integer $k$, real $\{a_j\}$ and complex $\{\nu_j\}$} subject to the conditions
\begin{itemize}
    \item $k\in\mathbb{N}_0\equiv\{0,1,2,...\}$
    \item \new{$\sum_{j=1}^N\left|a_j\right|\leq 2\pi$}
    \item $\sum_{j=1}^N\mathrm{Re}(\nu_j)>2k-N/2$ \item $\sum_{j=1}^N\mathrm{Re}(\nu_j)>2k-N/2+1$ \textit{if}:
    \begin{itemize}
        \item \new{$\sum_{j=1}^N\left|a_j\right|= 2\pi$} \textit{and/or}:
        \item There exists some vector \new{$\{s_j\}=\{\pm 1, \pm 1, ...\}$} such that $\sum_{j=1}^N s_ja_j=0$.
    \end{itemize}
\end{itemize}
We additionally assume $a_j\neq 0\,\forall\,j$, to avoid trivial results. \new{The range of validity is somewhat increased if some of $\{\nu_j\}$ are negative integers, in particular those with $j\in J$, for some non-empty set $J$. In this case,} \eqref{eq: main-result} applies also for integer $k\geq -\sum_{j'\in J}|\nu_{j'}|$.

\new{In practice, the condition that $\sum_{j=1}^N|a_j|\leq 2\pi$ is not a limitation. Assuming $\sum_{j=1}^N|a_j|= 2\pi A$ with $A\geq 1$, we may rescale $t\to \tilde t \equiv At$, $a_j\to \tilde a_j \equiv a_j/A$, giving
\beq
    \int_{0}^\infty dt\, t^{2k}\prod_{j=1}^N\left[t^{-\nu_j}J_{\nu_j}(a_jt)\right] &=& A^{\sum_j \nu_j-1-2k} \int_{0}^\infty d\tilde t\, \tilde t^{2k}\prod_{j=1}^N\left[\tilde t^{-\nu_j}J_{\nu_j}(\tilde a_j\tilde t)\right]\\\nonumber
    &=&A^{\sum_j \nu_j-1-2k}\sum_{m=0}^\infty \varepsilon_m\, m^{2k}\prod_{j=1}^N\left[m^{-\nu_j}J_{\nu_j}(a_jm/A)\right]
\eeq
noting that $\sum_{j=1}^N|\tilde a_j| = 2\pi A^{1-N}\leq 2\pi$.}

In the case $N=2$, labelling $\{a_j\} = \{a,b\}$, $\{\nu_j\} = \{\mu,\nu\}$, for real $a,b$, \eqref{eq: main-result} can be written
\beq\label{eq: rescaled-result}
    \int_{0}^\infty dt\, \frac{J_\mu(at)J_\nu(bt)}{t^{\mu+\nu-2k}} = \sum_{m=0}^\infty \varepsilon_m\,\frac{J_\mu(am)J_\nu(bm)}{m^{\mu+\nu-2k}}
\eeq
for $k\in\mathbb{N}_0$ (assuming neither $\mu$ nor $\nu$ are negative integers), $\mathrm{Re}(\mu+\nu)>2k-1$, and $|a|+|b|\leq 2\pi$. Additionally we require $\mathrm{Re}(\mu+\nu)>2k$ if $|a|+|b|=2\pi$ or $a=\pm b$. This matches the result of \citep{Dominici_2012}, but \new{with a larger domain of validity, for example including $a\geq\pi$.}

\section{Abel-Plana Theorem}\label{sec: abel-plana}
We now introduce the Abel-Plana theorem  \citep{abel_1823,plana_1820}, which will be used to prove our main result \eqref{eq: main-result} in \S\ref{sec: proof}. In the notation of \citep{ButzerPL2011TSFo}, the theorem states that
\beq\label{eq: abel-plana}
    \sum_{k=0}^\infty f(k) = \int_0^\infty dx\,f(x)+\frac{1}{2}f(0)+i\int_0^\infty dy\,\frac{f(iy)-f(-iy)}{e^{2\pi y}-1}
\eeq
where the function $f:\mathbb{C}\to\mathbb{C}$ obeys the following conditions:
\begin{enumerate}
    \item $f(z)$ is analytic in the closed half-plane $U = \{z\in\mathbb{C}:\mathrm{Re}(z)\geq 0\}$.
    \item $\lim_{y\to\infty}\left|f(x\pm iy)\right|e^{-2\pi y}=0$ uniformly in $x$ on every finite interval.
    \item $\int_0^\infty dy\,\left|f(x+ iy)-f(x-iy)\right|e^{-2\pi y}$ exists for every $x\geq 0$ and tends to zero as $x\to\infty$.
\end{enumerate}
The additional conditions
\begin{enumerate}[resume]
    \item $\int_0^\infty dx\,f(x)$ exists\new{,}
    \item $\lim_{n\to\infty}f(n)=0$\new{,}
\end{enumerate}
are often imposed, though \citep{ButzerPL2011TSFo} \new{consider} the second to be superfluous. The theorem itself is straightforwardly proved from the argument principle and Cauchy's integral theorem; \new{in essence, one considers the integral of $f(z)/[e^{-2\pi i z}-1]$, which has poles at integer $z$.} See \citep{olver1997asymptotics} and \citep{saharian2007generalized} for additional discussion of this theorem and its applications.


\section{Validity Conditions}\label{sec: conditions}
Consider the function $f:\mathbb{C}\to\mathbb{C}$:
\beq\label{eq: f-z}
    f(z) &=& z^{-\lambda}\prod_{j=1}^N J_{\nu_j}(a_jz)
\eeq
\new{where $J_\mu(z)$ is} a Bessel function of the first kind and $\mu_1,\mu_2,...,\mu_N$, $a_1,a_2,...,a_N$ and $\lambda$ \new{are complex parameters}. In the below, we demonstrate that $f(z)$ satisfies the Abel-Plana validity criteria given in \S\ref{sec: abel-plana}, subject to certain restrictions on each parameter. We will ignore the trivial cases in which at least one \new{element} of $\{a_j\}$ is zero.

\subsection{Condition 1}\label{subsec: cond-1}
\begin{quote}
    \new{\textit{$f(z)$ is analytic in the closed half-plane $U = \{z\in \mathbb{C}:\mathrm{Re}(z)\geq 0\}$.}}
\end{quote}
Since the Bessel functions are holomorphic \new{on} $\mathbb{C}$ except for a branch cut on the negative real axis, and $z^{-\lambda}$ is holomorphic in $\mathbb{C}\backslash\{0\}$, it follows that $f(z)$ is analytic in $U\backslash\{0\}$. \new{Near} $z=0$, we consider the asymptotic form for Bessel functions \new{with} $0<|z|\ll\sqrt{\nu+1}$:
\beq
    J_\nu(z) \approx \frac{1}{\Gamma(\nu+1)}\left(\frac{z}{2}\right)^\nu \quad \left(\nu \notin \mathbb{Z}^{-}\right)
\eeq
\citep[Eq.\,9.1.7]{abramowitz+stegun}, \new{where $\mathbb{Z}^{-}$ is the set of negative integers.} and $\Gamma(t)$ is the Gamma function. Thus, $f(z)$ has the asymptotic limit
\beq
    f(z) \approx \prod_{j=1}^N\left[\frac{1}{\Gamma(\nu_j+1)}\left(\frac{a_j}{2}\right)^{\nu_j}\right]\times\,z^{\sum_j \nu_j-\lambda}
\eeq
assuming no element of $\{\nu_j\}$ \new{is a} negative integer. The limit $z\to0$ (\textit{i.e.} $f(0)$) thus exists if $\sum_j \mathrm{Re}(\nu_j)-\mathrm{Re}(\lambda)>0$. \new{The limit} additionally exists if $\sum_j \nu_j-\lambda=0$. 

For negative integer order $\nu$, the Bessel function is instead approximated by
\beq
    J_\nu(z) \approx \frac{(-1)^\nu}{|\nu|!}\left(\frac{z}{2}\right)^{|\nu|} \quad \left(\nu \in \mathbb{Z}^{-}\right),
\eeq
\new{thus $\lim_{z\to0} J_\nu (z) = 0$ for all negative integer $\nu$. Given this,} the most general constraint on the existence of $f(0)$ is that
\beq\label{eq: f0-conditions}
    &&\sum_{j'\in J}|\nu_{j'}|+\sum_{j\notin J}\mathrm{Re}(\nu_j) - \mathrm{Re}(\lambda)>0\\\nonumber
    \text{or} && 
    \sum_{j'\in J}|\nu_{j'}|+\sum_{j\notin J}\nu_j - \lambda=0,
\eeq
\new{where $J$ is the set of indices} corresponding to negative integer $\nu_j$, \textit{i.e.} $J=\{j:\nu_{j}\in\mathbb{Z}^{-}\}$.
\new{Assuming \eqref{eq: f0-conditions} holds,} $f(z)$ is analytic on $U$ and its boundary, thus the condition is satisfied.

\subsection{Condition 2}\label{subsec: cond-2}
\begin{quote}
    \new{\textit{$\lim_{y\to\infty}|f(x\pm iy)|e^{-2\pi y}=0$ uniformly in $x$ on every finite interval.}}
\end{quote}

At large $|z|$, Bessel functions of the first kind have the following asymptotic form:
\beq\label{eq: bessel-asymptotic}
    J_\nu(z)\approx \sqrt{\frac{2}{\pi z}} \left[\cos\left(z-\frac{\nu \pi}{2}-\frac{\pi}{4}\right)+e^{\left|\mathrm{Im}(z)\right|}\mathcal{O}\left(|z|^{-1}\right)\right]
\eeq
\citep[Eq.\,9.2.1]{abramowitz+stegun}, assuming $|\mathrm{arg}(z)|<\pi/2$ (which is valid for all $z$ within \new{the closed half-plane} $U$). \new{To evaluate the large-$y$ limit of $|f(x\pm iy)|$, we require the term
\beq
    \left|J_\mu(az)\right| &\approx& \sqrt{\frac{1}{2\pi|az|}}\left|e^{i(az-\nu\pi/2-\pi/4)}-e^{-i(az-\nu\pi/2-\pi/4)}\right|,
\eeq
}ignoring the subdominant $\mathcal{O}\left(|z|^{-1}\right)e^{|\mathrm{Im}(z)|}$ term (appropriate given that $|\cos z|$ is $\mathcal{O}(e^{|\mathrm{Im}(z)|})$ for large $|\mathrm{Im}(z)|$). Writing $z = x\pm iy$, \new{the dominant term is the one involving $y$ with a positive coefficient;
\beq
    \left|J_\mu(a[x\pm iy])\right| \approx \sqrt{\frac{1}{2\pi|a||y|}}\mathrm{exp}\left(\mathrm{sgn}[\mathrm{Re}(a)]\left[\mathrm{Re}(a)y\pm \frac{\pi}{2}\mathrm{Im}(\nu)\mp\mathrm{Im}(a)x\right]\right),
\eeq
noting that $|e^{i\phi}|=1$ for all real $\phi$. Utilizing definition \eqref{eq: f-z},} we find the following asymptotic limit of $|f(x\pm iy)|$ for large $y>0$:
\beq\label{eq: abs-fx-asmypt}
    \left|f(x\pm iy)\right|&\approx&y^{-[N/2+\mathrm{Re}(\lambda)]}\exp\left[\pm\frac{\pi}{2}\mathrm{Im}(\lambda)\right]\\\nonumber
        &&\,\quad\,\times\,\prod_{j=1}^N\left[(2\pi|a_j|)^{-1/2}\exp\left(\mathrm{sgn}[\mathrm{Re}(a_j)]\left[\mathrm{Re}(a_j)y\pm \frac{\pi}{2}\mathrm{Im}(\nu_j)\mp\mathrm{Im}(a_j)x\right]\right)\right],
\eeq
noting that \neww{$|z^\lambda| = |z|^{\mathrm{Re}(\lambda)}\mathrm{exp}\left[-\mathrm{Im}(\lambda)\mathrm{arg}(z)\right]$ with $\mathrm{arg}(z)\approx\pm\pi/2$.}


We may now test the condition:
\beq\label{eq: ab-f-result}
    \lim_{y\to\infty}\left|f(x\pm iy)\right|e^{-2\pi y} &=& \exp\left[\pm\frac{\pi}{2}\mathrm{Im}(\lambda)\right]\prod_{j=1}^N\left\{(2\pi|a_j|)^{-1/2}\exp\left[\pm\frac{\pi}{2}\mathrm{sgn}[\mathrm{Re}(a_j)]\mathrm{Im}(\nu_j)\right]\right\}\nonumber\\
    &&\,\times\,\exp\left[\mp x\sum_{j=1}^N\mathrm{sgn}[\mathrm{Re}(a_j)]\mathrm{Im}(a_j)\right]\\\nonumber
    &&\,\times\,\lim_{y\to\infty}\left[y^{-[N/2+\mathrm{Re}(\lambda)]}\exp\left(y\left[\sum_{j=1}^N|\mathrm{Re}(a_j)|-2\pi\right]\right)\right].
\eeq
The limit is zero if (i) $\sum_{j=1}^N|\mathrm{Re}(a_j)|<2\pi$ or (ii) $\sum_{j=1}^N|\mathrm{Re}(a_j)|=2\pi$ and $\mathrm{Re}(\lambda)>-N/2$. Furthermore, if $\sum_{j=1}^N\mathrm{sgn}[\mathrm{Re}(a_j)]\mathrm{Im}(a_j)=0$, the limit is uniformly approached in $x$, as required. Henceforth, we will assume the stronger condition that all $a_j$ are real; this is later required in \S\ref{sec: proof}.

\subsection{Condition 3}\label{subsec: cond-3}
\begin{quote}
    \new{\textit{$\int_0^\infty dy\,|f(x + iy)-f(x-iy)|e^{-2\pi y}$ exists for every $x\geq 0$ and tends to zero as $x\to\infty$.}}
\end{quote}
For this, we first note that the quantity of interest is bounded on both sides: 
\beq\label{eq: condition-3-tmp}
    0< \int_0^\infty dy\,|f(x+iy)-f(x-iy)|e^{-2\pi y}\leq \int_0^\infty dy\,|f(x+iy)|e^{-2\pi y}+\int_0^\infty dy\,|f(x-iy)|e^{-2\pi y}
\eeq
using the triangle inequality and noting that $|f(x+iy)-f(x-iy)|\geq 0$ and \new{$e^{-2\pi y}> 0$} for all $y\geq 0$. From \new{condition 1 (\S\ref{subsec: cond-1})}, $f(x\pm iy)$ is holomorphic for all $x\geq 0$, and thus Riemann-integrable across any finite interval $y\in[u,v]$ with $u,v\geq 0$. To show that the infinite integral exists, we must consider its behavior at infinity. For sufficiently large $y$, $|f(x\pm iy)|$ can be replaced by \eqref{eq: ab-f-result}, such that
\beq
    |f(x\pm iy)|e^{-2\pi y} \approx K_0 y^{-[N/2+\mathrm{Re}(\lambda)]}\exp\left(-\left[2\pi-\sum_{j=1}^N |a_j|\right]y\right)
\eeq
for some $K_0>0$, assuming \new{that $\{a_j\}$ are real.} If $\sum_{j=1}^N|a_j|<2\pi$, then there exist real numbers $c_1, K_1$ such that
\beq
    |f(x\pm iy)|e^{-2\pi y}< K_1e^{-c_1y}
\eeq
for all $y>M0$ with \new{$M\gg 0$}, $0<c_1<2\pi-\sum_{j=1}^{N}|a_j|$ and $K_1>0$. This implies
\beq
    \int_M^\infty dy\,|f(x\pm iy)|e^{-2\pi y}<K_1\int_M^\infty dy\,e^{-c_1y} = \frac{K_1}{c_1}e^{-c_1M}<\infty.
\eeq
If instead we assume $\sum_{j=1}^N |a_j|=2\pi$ and $\mathrm{Re}(\lambda)>-N/2$, there must exist real numbers $c_2,K_2$ such that
\beq
    |f(x\pm iy)|e^{-2\pi y}< K_2 e^{-c_2}
\eeq
for all $y>M\gg 0$, with $0<c_2<\mathrm{Re}(\lambda)+N/2$. In this case, 
\beq
    \int_M^\infty dy\,|f(x\pm iy)|e^{-2\pi y}<K_2\int_M^\infty dy\,y^{-c_2} = \frac{K_2}{1-c_2}M^{1-c_2}<\infty
\eeq
if $1-c_2<0$, and thus $\mathrm{Re}(\lambda)>1-N/2$, which is a slightly stronger bound than before. 

If such restrictions are satisfied, the integral of $|f(x\pm iy)|e^{-2\pi y}$ over the whole range $y\in[0,\infty]$ must be finite, and thus, by \eqref{eq: condition-3-tmp}, $\int_0^\infty dy\,|f(x+iy)-f(x-iy)|e^{-2\pi y}$ must exist. Furthermore, \eqref{eq: abs-fx-asmypt} demonstrates that $|f(x\pm iy)|e^{-2\pi y}$ is suppressed by $x^{-(N/2+\mathrm{Re(\lambda)})}$ at large $x$; thus, if $\mathrm{Re}(\lambda)>-N/2$, the integral tends to zero as $x\to\infty$ \new{and the condition is satisfied.}

\subsection{Condition 4}\label{subsec: cond-4}
\begin{quote}
    \new{\textit{$\int_0^\infty dx\,f(x)$ exists.}}
\end{quote}
First we note that 
\beq
    \left|\int_0^\infty dx\,f(x)\right|\leq \int_0^\infty dx\,|f(x)|;
\eeq
thus proving the latter to exist is a sufficient (but not necessary) condition for the existence of $\int_0^\infty dx\,f(x)$. Proof of this proceeds analogously to \new{condition 3 (\S\ref{subsec: cond-3})}, where we first note that, assuming \new{condition 1 (\S\ref{subsec: cond-1})} be satisfied, $f(x)$ is Riemann-integrable on any finite interval. Secondly, using the asymptotic form of $|f(x)|$ from the first line of \eqref{eq: abs-fx-asmypt}, we note that there must exist constants $K_3$, $c_3$ such that
\beq
    |f(x)|< K_3x^{-c_3}
\eeq
for all $x>M'\gg 0$, where $0<c_3<\mathrm{Re}(\lambda)+N/2$. Using this
\beq
    \int_{M'}^\infty dx\,|f(x)|<K_3\int_{M'}^\infty dx\,x^{-c_3} = \frac{K_3}{1-c_3}(M')^{1-c_3}<\infty
\eeq
if $1-c_3<0$ and thus $\mathrm{Re}(\lambda)>1-N/2$.

In practice, a stronger bound is in fact possible in most cases, \new{since the integrand is oscillatory}. To show this, we again consider $x>M'$ for large positive $M'$, with the asymptotic form
\beq\label{eq: fx-asympt}
    f(x)\approx x^{-\lambda-N/2}\prod_{j=1}^N \left[(2\pi a_j)^{-1/2}\cos\left(a_jx-\nu_j\frac{\pi}{2}-\frac{\pi}{4}\right)\right]
\eeq
from \eqref{eq: bessel-asymptotic}, \new{with $\mathrm{Im}(z) =0$}. Next, note that the product of two \new{cosines} can be rewritten as a sum via the standard relation
\beq
    \cos\theta \cos\phi = \frac{1}{2}\left[\cos(\theta+\phi)+\cos(\theta-\phi)\right]
\eeq
for arbitrary complex $\theta$, $\phi$. In this way, we may reduce the product of $N$ cosine functions into a \new{linear sum of} $2^N$ terms of the form $\cos(Ax+B)$. In particular the $A$ coefficients \new{are of the form} $A = \sum_{j=1}^N s_j a_j$ where $s_j = \{\pm1, \pm1, ...\}$. Each term is thus an oscillatory function of $x$, provided that \new{one} cannot find an vector $\{s_j\}$ satisfying $\sum_j s_j a_j=0$. If one \textit{can} be found, then we obtain a non-oscillatory term. Physically, this corresponds to a case when one of the `beat-frequencies' is equal to zero. \new{In this case, $\int_0^\infty dx\,f(x)$ contains a logarithmic divergence, and the condition is not satisfied.}

Taking the above approach, we can write schematically
\beq\label{eq: fx-asympt2}
    f(x)&\approx& K_4\,x^{-\lambda-N/2} \sum_{i=1}^{2^N}\cos\left(A_ix+B_i\right)\\\nonumber
    \Rightarrow \int_{M'}^\infty dx\,f(x) &\approx & K_4\sum_{i=1}^{2^N} \int_{M'}^\infty dx\,x^{-\lambda-N/2}\cos\left(A_ix+B_i\right)
\eeq
for some $K_4$. Assuming \new{each} $A_i\neq 0$, each integral may be divided up into regions where the integrand is negative and positive:
\beq
    \int_{M^*}^\infty dx\,x^{-\lambda-N/2}\cos\left(A_ix+B_i\right) &=& \sum_{p=0}^\infty \int_{M^*+p\pi/A_i}^{M^*+(p+1)\pi/A_i} dx\,x^{-\lambda-N/2}\cos\left(A_ix+B_i\right)\\\nonumber
    &=&\sum_{p=0}^\infty (-1)^p \int_{M^*+p\pi/A_i}^{M^*+(p+1)\pi/A_i} dx\,\left|x^{-\lambda-N/2}\cos\left(A_ix+B_i\right)\right|
\eeq
where $M^*$ is the first increasing zero of $\cos(A_ix + B_i)$ with $M^*>M'$. Provided $\mathrm{Re}(\lambda)>-N/2$, the magnitude of the integrand decreases monotonically with $p$ (and to zero as $p\to\infty$) thus the alternating series, and hence the integral, converges by the Leibniz criterion. Using this approach, we find that $\int_0^\infty dx\,f(x)$ exists for all $\mathrm{Re}(\lambda)>-N/2$, provided that we exclude $A_i=0$, \textit{i.e.} provided that there does \new{not} exist $\{s_j\}=\{\pm,\pm,...\}$ such that $\sum_{j=1}^N s_j a_j=0$. If \new{such a vector \textit{can} be found}, we \new{instead} require $\mathrm{Re}(\lambda)>1-N/2$ \new{to avoid logarithmic divergences in the integral}.

Since \new{the above requirements are} somewhat non-trivial, it is useful to \new{consider} the simpler case $N=2$, \new{which has} $\{a_j\} = \{a,b\}$, $\{\nu_j\} = \{\mu,\nu\}$ \new{and was discussed in \citep{Dominici_2012}}. Here, $\int_0^\infty dx\,f(x)$ is a known integral that can be written in terms of hypergeometric functions \citep[e.g.\,Eq.\,6.574,][]{1994tisp.book.....G}. \new{That solution} requires $\mathrm{Re}(\lambda)>-1$ and real $a>0$, $b>0$, with the additional constraint $\mathrm{Re}(\lambda)>0$ if $a=b$. Given that $J_\nu(-ax)=(-1)^\nu J_\nu(ax)$ for real $ax$ (as here), the \new{regime of validity} may be extended also to $a<0$ and $b<0$, and is \new{subject} to the same limits as our approach.

\subsection{Condition 5}\label{subsec: cond-5}
\begin{quote}
    \new{\textit{$\lim_{n\to\infty}f(n) = 0$.}}
\end{quote}
From \eqref{eq: fx-asympt}, the behavior of \new{$f(n)$} at large \new{(integer)} $n$ is that of a cosine bounded by the envelope $n^{-[\mathrm{Re}(\lambda)+N/2]}$. For real $\{a_j\}$, the magnitude of each cosine is bounded by a constant (equal to unity if $\nu_j$ is real), thus $|f(n)|\leq K_5\,n^{-[\mathrm{Re}(\lambda)+N/2]}$ for some constant $K_5>0$. This implies
\beq
    0\leq \left|\lim_{n\to\infty} f(n)\right| \leq K_5\lim_{n\to\infty} n^{-[\mathrm{Re}(\lambda)+N/2]}.
\eeq
Assuming $\mathrm{Re}(\lambda)>-N/2$, the limit is zero, satisfying the condition.


\section{Proof of the Integral-to-Sum Relation}\label{sec: proof}
We now insert our definition of $f(z)$ \eqref{eq: f-z} into the Abel-Plana theorem \eqref{eq: abel-plana}. This gives
\beq
    \sum_{m=0}^\infty  m^{-\lambda}\prod_{j=1}^NJ_{\nu_j}(a_jm) &=& \int_{0}^\infty dt\, t^{-\lambda}\prod_{j=1}^NJ_{\nu_j}(a_jt) + \frac{1}{2}\lim_{s\to0} s^{-\lambda}\prod_{j=1}^NJ_{\nu_j}(a_js)\\\nonumber
    &&\,+\,i\int_0^\infty \frac{dy}{e^{2\pi y}-1}\left[(iy)^{-\lambda}\prod_{j=1}^NJ_{\nu_j}(ia_jy)-(-iy)^{-\lambda}\prod_{j=1}^NJ_{\nu_j}(-ia_jy)\right]
\eeq
By introducing the $\varepsilon_m$ factor of \eqref{eq: epsilon-m-def}, we may shift the second term on the RHS of \new{the first line} to the LHS. \new{Secondly, we note that $J_\mu(it)\equiv i^\mu I_\mu(t)$, where $I_\mu$ is a modified Bessel function of the first kind. This leads} to
\beq
    \sum_{m=0}^\infty \varepsilon_m\,m^{-\lambda}\prod_{j=1}^NJ_{\nu_j}(a_jm) &=& \int_{0}^\infty dt\, t^{-\lambda}\prod_{j=1}^NJ_{\nu_j}(a_jt) \\\nonumber
    &&\,+\,i\int_0^\infty dy\,\frac{(iy)^{-\lambda}}{e^{2\pi y}-1}\left[\prod_{j=1}^NI_{\nu_j}(a_jy)-(-1)^{-\lambda}\prod_{j=1}^NI_{\nu_j}(-a_jy)\right].
\eeq
\new{For} real arguments $t$ (as here, since \new{we} assume all $a_j$ to be real), $I_\nu(t) \equiv (-t)^\nu t^{-\nu}I_\nu(t)\equiv (-1)^\nu I_\nu(t)$, thus
\beq
    \sum_{m=0}^\infty \varepsilon_m\,m^{-\lambda}\prod_{j=1}^NJ_{\nu_j}(a_jm) &=& \int_{0}^\infty dt\, t^{-\lambda}\prod_{j=1}^NJ_{\nu_j}(a_jt) \\\nonumber
    &&\,+\,i\int_0^\infty dy\,\frac{(iy)^{-\lambda}}{e^{2\pi y}-1}\left[1-(-1)^{-\lambda+\sum_{j=1}^N\nu_j}\right]\prod_{j=1}^NI_{\nu_j}(a_jy).
\eeq
If $\sum_{j=1}^N \nu_j -\lambda = 2k$ for $k\in\mathbb{Z}$, the second line vanishes identically and we obtain the desired result \eqref{eq: main-result}. 

We finally consider the \new{restrictions on $\lambda$, $\{a_j\}$ and $\{\nu_j\}$ arising from} the Abel-Plana validity conditions of \S\ref{sec: abel-plana}. If none of $\nu_j$ are negative integers, condition 1 (\S\ref{subsec: cond-1}) necessitates $\mathrm{Re}\left(\sum_{j=1}^N\nu_j-\lambda\right)>0$ unless $\sum_{j=1}^N\nu_j-\lambda=0$; here, $\sum_{j=1}^N\nu_j-\lambda=2k$, thus the condition is satisfied for all $k\in\mathbb{N}_0$. If some of $\{\nu_j\}$ \textit{are} negative integers, in \new{those with $j$ in some set $J$}, we instead find a slightly more lenient condition
\beq
    \sum_{j'\in J}(|\nu_{j'}|-\nu_{j'})+2k\geq 0 \quad \Rightarrow \quad k\geq -\sum_{j'\in J}|\nu_{j'}|.
\eeq
From \new{condition 3 (\S\ref{subsec: cond-3})}, $\mathrm{Re}(\lambda)>-N/2$ and thus $\sum_{j=1}^N\mathrm{Re}(\nu_j)-2k>-N/2$. Furthermore, we require $\sum_{j=1}^N|a_j|< 2\pi$ unless $\sum_{j=}^N|a_j| = 2\pi$ and $\mathrm{Re}(\lambda)>1-N/2$: this implies that $\sum_{j=1}^N\mathrm{Re}(\nu_j)-2k>1-N/2$. Finally, from \new{condition 4 (\S\ref{subsec: cond-4})}, $\sum_{j=1}^N\mathrm{Re}(\nu_j)-2k>1-N/2$ is also required if the exists a vector $\{s_j\} = \{\pm 1, \pm 1, ...\}$ such that $\sum_{j=1}^N s_j a_j=0$. These assumptions additionally satisfy the remaining conditions, thus we arrive at the result of \S\ref{sec: main-result}.

\section{Alternate Proof via the Poisson Summation Formula}\label{sec: alternate-proof}
\neww{Below, we sketch an alternative proof of \eqref{eq: main-result}, using the Poisson summation formula.\footnote{We thank Jeremy Goodman for suggesting this approach.} For this, first consider a general function $g:\mathbb{R}\to\mathbb{C}$. Denoting the Fourier transform of $g$ by $\hat{g}$, the Poisson summation formula links the infinite sum of $g$ and $\hat{g}$:
\beq
    h\sum_{n=-\infty}^\infty g(nh) = \sum_{k=-\infty}^\infty \hat{g}(k/h)
\eeq
for arbitrary $h\in\mathbb{R}$ \citep[e.g.,][]{zygmund}, assuming that $\hat{g}$ exists. If $g$ is band-limited, such that $\hat{g}(p)$ has support only for $p\in(-h,h)$, then the RHS reduces to $\hat{g}(0)$. Inserting the integral definition of $\hat{g}(0)$, this gives
\beq\label{eq: poisson2}
    h\sum_{n=-\infty}^\infty g(nh) = \hat{g}(0) = \int_{-\infty}^\infty dt\,g(t).
\eeq
}

\neww{To apply this result in our context, we must first ascertain whether $f(t)$ \eqref{eq: f-z} is band-limited. For this purpose, we first note that the Bessel functions may be written in integral form as
\beq
    J_\nu(t) = \frac{\left(t/2\right)^\nu}{\sqrt{\pi}\,\Gamma\left(\nu+\tfrac{1}{2}\right)}\int_0^\pi d\theta\,\left(\sin\theta\right)^{2\nu}\cos\left(t\cos\theta\right)
\eeq
for $\mathrm{Re}(\nu)>-1/2$ \citep[\S10.9.4]{nist_dlmf}. Since the RHS is a weighted average of cosines in $t$ with frequencies in the range $[-1/(2\pi),1/(2\pi)]$, it follows that the Fourier transform of $J_\nu(t)$ has support only over this range. The product of $N$ Bessel functions with scaling parameters $\{a_j\}$ can thus be written in a form involving the product of $N$ cosines in $t$, each of which have frequencies in the range $[-a_j/(2\pi),a_j/(2\pi)]$. Since $t^{-\lambda}$ is a polynomial, and thus of zero bandwidth, the function $f(z)$ \eqref{eq: f-z} contains only frequencies in the range $[-a_j/(2\pi),a_j/(2\pi)]$ and will be band-limited by $[-1,1]$, provided that $\sum_j a_j < 2\pi$ (for real $a_j$).\footnote{This is true also for $\sum_j a_j = 2\pi$ if there are no beat-frequencies, as in \S\ref{subsec: cond-4}.} Applying the Poisson summation formula with $h=1$ leads to
\beq
    \sum_{m=-\infty}^\infty f(m) = \hat{f}(0) = \int_{-\infty}^\infty dt\,f(t).
\eeq
We can remove the $m<0$ and $t<0$ terms by relabelling:
\beq\label{eq: poisson2}
    \frac{1}{2}\sum_{m=0}^\infty \varepsilon_m\left[f(m)+f(-m)\right] = \frac{1}{2}\hat{f}(0) = \frac{1}{2}\int_{0}^\infty dt\,\left[f(t)+f(-t)\right],
\eeq
introducing the $\varepsilon_m$ coefficient \eqref{eq: epsilon-m-def} to capture the special case $m=0$. Finally, we insert the definition of $f$, giving
\beq\label{eq: poisson-result}
    \sum_{m=0}^\infty \varepsilon_m m^{-\lambda}\left[1+(-1)^{-\lambda-\sum_j=1^N\nu_j}\right] \prod_{i=1}^N\left[J_{\nu}(a_jm)\right]= \int_{0}^\infty dt\,t^{-\lambda}\left[1+(-1)^{-\lambda+\sum_j \nu_j}\right]\prod_{j=1}^N\left[J_{\nu}(a_jt)\right],
\eeq
where we have noted that $J_{\nu}(-at) = (-1)^\nu J_{\nu}(at)$ \citep[\S10.11.1]{nist_dlmf}. Setting $\lambda = \sum_{j=1}^N\nu_j - 2k$ for $k\in\mathbb{N}_0$ as before, we obtain the desired result \eqref{eq: main-result}. The remaining conditions in \S\ref{sec: conditions} ensure that (a) the Poisson summation formula is valid, and (b) the integral and sum in \eqref{eq: poisson-result} exist.
}

\section{Summation Convergence}
\neww{We briefly comment on the convergence of the infinite sum appearing in \eqref{eq: main-result}. This is important for assessing the utility of the result as a method to evaluate Bessel function integrals. For this purpose, we consider truncating $\sum_{m=0}^\infty f(m)$ at $m=M$ for large $M$. By the integral test for convergence, $\sum_{m=M}^\infty f(m)$ converges absolutely \textit{iff} $\int_{M}^\infty dx\,|f(x)|$ converges. Condition 5 (\S\ref{subsec: cond-5}) implies that
\beq
    \int_{M}^\infty dx\,|f(x)| \leq K\int_{M}^\infty dx\, x^{-[\mathrm{Re}(\lambda)+N/2]}.
\eeq
The integral exists for all $\mathrm{Re}(\lambda)+N/2>1$, \textit{i.e.} the sum is absolutely convergent for $\sum_j \mathrm{Re}(\nu_j)>2k-N/2+1$. Furthermore, the integral test also states that
\beq
    \sum_{m=M}^\infty |f(m)| &\leq& f(M) + \int_{M}^\infty dx\,f(x)= KM^{-[\mathrm{Re}(\lambda)+N/2]}\left[1 + \frac{M}{1-[\mathrm{Re}(\lambda)+N/2]}\right]\\\nonumber 
    &=& \mathcal{O}(M^{1-[\mathrm{Re}(\lambda)+N/2]}),
\eeq
bounding the truncation error.}

\neww{Considering the case $0<\mathrm{Re}(\lambda)+N/2<1$, the sum $\sum_{m=0}^\infty f(m)$ is instead \textit{conditionally} convergent. This can be shown in a manner analogous to the derivation of condition 4 (\S\ref{subsec: cond-4}), first separating out the various Fourier frequencies (cf.\,\ref{eq: fx-asympt2}):
\beq
    \sum_{m=M}^\infty f(m) = K' \sum_{m=M}^\infty m^{-\lambda -N/2}\sum_{i=1}^{2^N}\cos(A_im+B_i),
\eeq
for some $K'$ and sufficiently large $M$, where the validity conditions of \S\ref{sec: conditions} ensure that $A_i\neq 0$. Considering a single frequency $A_i$ (\textit{i.e.} one of the $2^N$ terms in the summation over $i$), the summation over $m$ may be recast as an alternating series:
\beq\label{eq: summation-tmp}
    K' \sum_{m=M_*}^\infty m^{-\lambda -N/2}\cos(A_im+B_i) = \sum_{p=0}^\infty (-1)^p \sum_{m=M_*+\left\lceil{p\pi/A_i}\right\rceil}^{M_*+\left\lfloor{(p+1)\pi/A_i}\right\rfloor} m^{-\lambda -N/2}\cos(A_im+B_i)
\eeq
where we start at $M_*>M$, the first increasing zero of $\cos(A_im+B_i)$ beyond $M$. Note that the real part of the summand is explicitly positive. Each term in the $p$ summation is positive and decreasing, and, assuming $\mathrm{Re}(\lambda)>-N/2$, \eqref{eq: summation-tmp} converges by the alternating series test. Truncating at $M_*$, the error in the term with frequency $A_i$ cannot be greater than the first excluded $p$ term, \textit{i.e.}
\beq
    \left|\sum_{m=M_*}^\infty m^{-\lambda -N/2}\cos(A_im+B_i)\right| \leq  \sum_{m=M_*}^{M_*+\left\lfloor{\pi/A_i}\right\rfloor} \left|m^{-\lambda -N/2}\cos(A_im+B_i)\right|\leq K''M_*^{-[\mathrm{Re}(\lambda)+N/2]}
\eeq
for some $K''>0$. This applies to all $2^N$ choices of frequency $A_i$; combining, we find
\beq
    \left|\sum_{m=M}^\infty f(m)\right| = \mathcal{O}\left(M^{-[\mathrm{Re}(\lambda)+N/2]}\right),
\eeq
giving the relevant truncation error.
}

\section{Practical Demonstration \& Conclusions}\label{sec: demonstration}

\begin{figure*}[t!]
    \centering
    \begin{subfigure}[t]{0.5\textwidth}
        \centering
        \includegraphics[width=\textwidth]{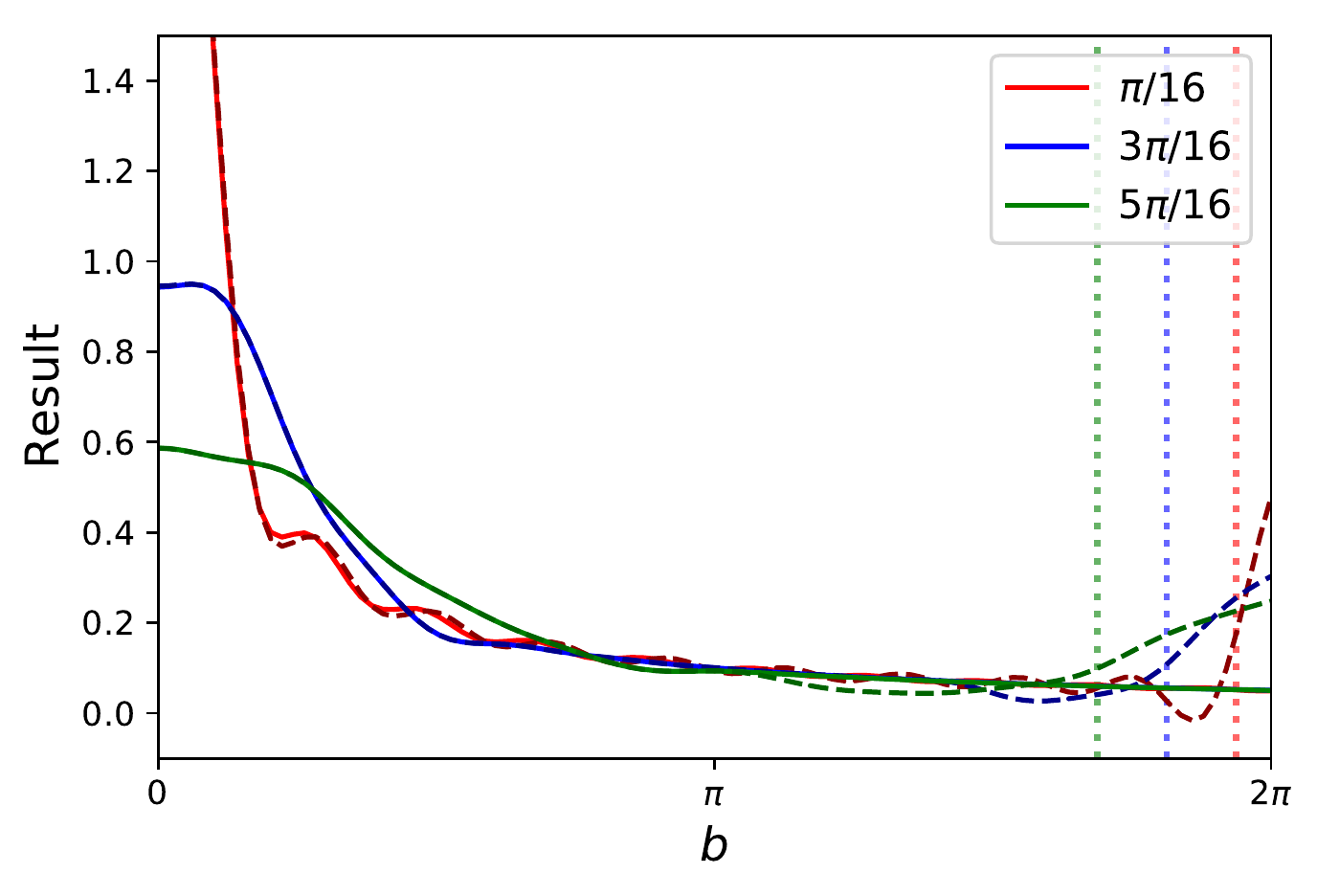}
        \caption{$\int_0^\infty dt\,t^{-2}J_{1/2}(at)J_{3/2}(bt)$}
    \end{subfigure}%
    ~ 
    \begin{subfigure}[t]{0.5\textwidth}
        \centering
        \includegraphics[width=\textwidth]{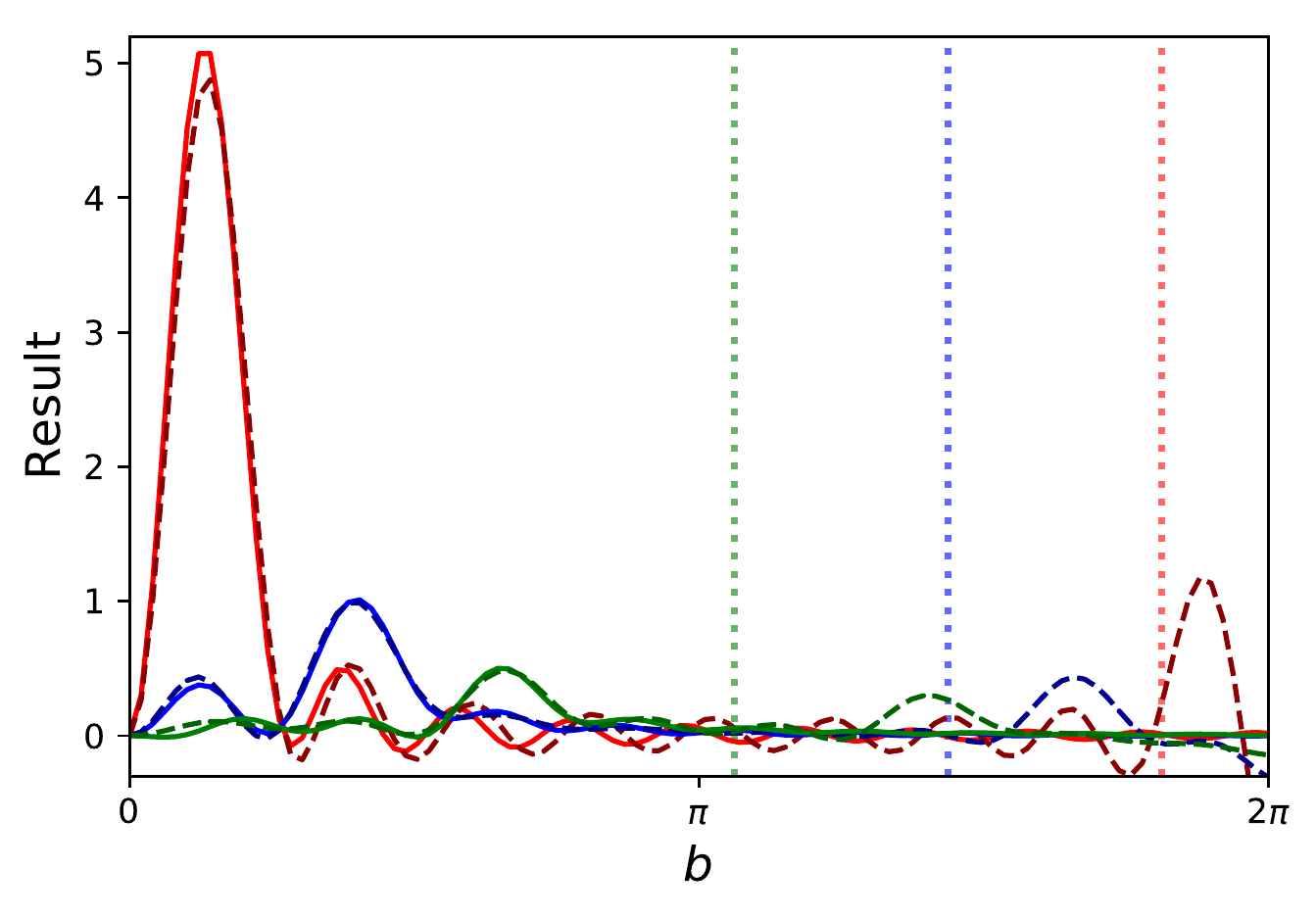}
        \caption{$\int_0^\infty dt\,t\,J_{0}(at)J_{1}(at)J_{2}(bt)$}
    \end{subfigure}\\
    \begin{subfigure}[t]{0.5\textwidth}
        \centering
        \includegraphics[width=\textwidth]{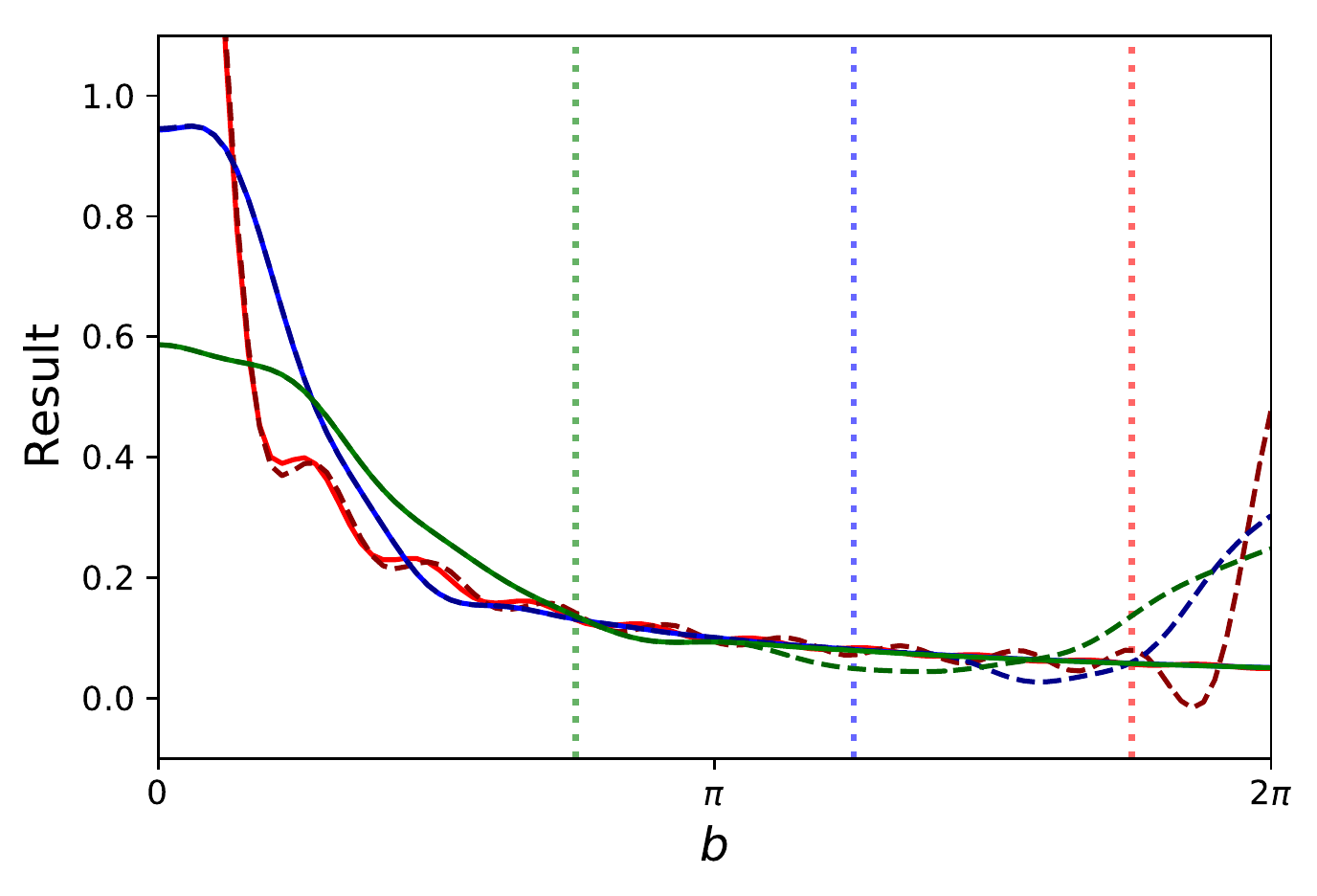}
        \caption{$\int_0^\infty dt\,J_{-3/2}(at)J_{-1}(at)J_{1/2}(at)J_0(bt)$}
    \end{subfigure}
    \caption{Comparison of Bessel function integrals computed via numerical quadrature (light solid lines) and the discrete summations of this work using ten terms (\ref{eq: main-result}, dark dashed lines). Results are shown for products of two, three and four Bessel functions, as detailed in the captions. The $a_j$ coefficients of $N-1$ Bessel function arguments are fixed to $a = \pi/16$ (red), $3\pi/16$ (blue) or $5\pi/16$ (green), but the $N$-th, denoted by $b$, is varied. The coefficients satisfy the validity criteria of \S\ref{sec: main-result} for all $b$ to the left of the dotted vertical lines, thus we expect the sums and integrals to agree in the limit of infinite sampling points. Note that case (b) has $\sum_j \mathrm{Re}(\nu_j) > 2k-N/2+1$, and is thus not expected to converge at $\sum_j a_j = 2\pi$; \neww{further, it is only conditionally convergent}. Furthermore, case (c) has $k=-1$; possible since one of the Bessel functions has a negative integer argument. We find good convergence in all cases; this can be improved further by including more points in the discrete summation or increasing the domain limit (set here to $t\in[0,10]$).}\label{fig: numerics}
\end{figure*}

Fig.\,\ref{fig: numerics} presents a demonstration of our main result \eqref{eq: main-result}. For this purpose, we consider three representative choices of hyperparameters $\{a_j\}$, $\{\nu_j\}$, $N$ and $k$, which satisfy the conditions of \S\ref{sec: main-result}. In each case, we compute both the LHS and RHS of \eqref{eq: main-result} separately, using numerical quadrature to perform the integral (truncating at $t_\mathrm{max} = 10$) and evaluate the sum including the first $10$ terms. In the regime $\sum_{j=1}^N |a_j|<2\pi$, the results are found to be in good agreement; discrepancies can be reduced by using a greater number of terms or a larger $s_\mathrm{max}$. This is additionally true for the special case of negative integer $\nu_j$ (since Fig.\,\ref{fig: numerics}c converges with $k=-1$), \neww{$\sum_j \mathrm{Re}(\nu_j) > 2k-N/2+1$, which is only conditionally convergent}. For $\sum |a_j|>2\pi$, the summation does not converge to the integral result, as expected. In practice, this can be avoided via \eqref{eq: rescaled-result} which extends the validity to all $\sum_{j=1}^N |a_j|$.

\vskip 4pt
In conclusion, we find that, subject to a number of assumptions, we can derive a useful formula relating the sums and integrals of products of $N$ Bessel functions of both real and complex order, allowing straightforward evaluation of highly oscillatory integrals. The result is straightforward to prove using the Abel-Plana theorem, and agrees with previous work for $N=2$, but with a significantly simpler proof. Whilst we consider only spherical Bessel functions in this work, we expect similar results to apply to other functions, provided they obey mild growth conditions and the relation $f(iy) = f(-iy)$ for real $y$.

\section*{Acknowledgements}
\small 
We thank \neww{Robert Cahn, Jeremy Goodman, Jiamin Hou, Sean Lake, Kiersten Meigs, David Spergel} and William Underwood for insightful discussions. OP acknowledges funding from the WFIRST program through NNG26PJ30C and NNN12AA01C.

\bibliographystyle{unsrt}
\bibliography{refs}

\begin{thebibliography}{10}

\bibitem{1993PhRvA..48.1134B}
P.~M. {Bergstrom}, T.~{Suri{\'c}}, K.~{Pisk}, and R.~H. {Pratt}.
\newblock {Compton scattering of photons from bound electrons: Full
  relativistic independent-particle-approximation calculations}.
\newblock {\em \pra}, 48(2):1134--1162, August 1993.

\bibitem{1995PhRvC..51.2031M}
R.~{Mehrem}, J.~T. {Londergan}, and G.~E. {Walker}.
\newblock {Isolating physical effects in the exclusive (N,N'{\ensuremath{\pi}})
  reaction}.
\newblock {\em \prc}, 51(4):2031--2043, April 1995.

\bibitem{2007PhRvD..76h3523F}
J.~R. {Fergusson} and E.~P.~S. {Shellard}.
\newblock {Primordial non-Gaussianity and the CMB bispectrum}.
\newblock {\em \prd}, 76(8):083523, October 2007.

\bibitem{2009PhRvD..79f3502H}
A.~{Hohenegger}.
\newblock {Solving the homogeneous Boltzmann equation with arbitrary scattering
  kernel}.
\newblock {\em \prd}, 79(6):063502, March 2009.

\bibitem{2010JOSAA..27.2144B}
Natalie {Baddour}.
\newblock {Operational and convolution properties of three-dimensional Fourier
  transforms in spherical polar coordinates}.
\newblock {\em Journal of the Optical Society of America A}, 27(10):2144,
  September 2010.

\bibitem{2012PhRvC..86a4003W}
K.~A. {Wendt}, R.~J. {Furnstahl}, and S.~{Ramanan}.
\newblock {Local projections of low-momentum potentials}.
\newblock {\em \prc}, 86(1):014003, July 2012.

\bibitem{2012CQGra..29v4007A}
D.~{Ahn}.
\newblock {Black hole state evolution, final state and Hawking radiation}.
\newblock {\em Classical and Quantum Gravity}, 29(22):224007, November 2012.

\bibitem{2013MNRAS.434.1808M}
Philipp~M. {Merkel} and Bj{\"o}rn~Malte {Sch{\"a}fer}.
\newblock {Intrinsic alignments and 3d weak gravitational lensing}.
\newblock {\em \mnras}, 434(2):1808--1820, September 2013.

\bibitem{2015LRR....18....1C}
Vitor {Cardoso}, Leonardo {Gualtieri}, Carlos A.~R. {Herdeiro}, and Ulrich
  {Sperhake}.
\newblock {Exploring New Physics Frontiers Through Numerical Relativity}.
\newblock {\em Living Reviews in Relativity}, 18(1):1, September 2015.

\bibitem{2015OExpr..2333044F}
Ivan {Fernandez-Corbaton}, Stefan {Nanz}, Rasoul {Alaee}, and Carsten
  {Rockstuhl}.
\newblock {Exact dipolar moments of a localized electric current distribution}.
\newblock {\em Optics Express}, 23(26):33044, December 2015.

\bibitem{2015PhRvD..91b3508V}
Zvonimir {Vlah}, Uro{\v{s}} {Seljak}, and Tobias {Baldauf}.
\newblock {Lagrangian perturbation theory at one loop order: Successes,
  failures, and improvements}.
\newblock {\em \prd}, 91(2):023508, January 2015.

\bibitem{2016JOSAB..33C..30S}
P.~{Springer}, S.~W. {Koch}, and M.~{Kira}.
\newblock {Excitonic terahertz absorption in semiconductors with effective-mass
  anisotropies}.
\newblock {\em Journal of the Optical Society of America B Optical Physics},
  33(7):C30, July 2016.

\bibitem{2016TMP...188.1197K}
A.~V. {Kisselev}.
\newblock {Approximate formulas for moderately small eikonal amplitudes}.
\newblock {\em Theoretical and Mathematical Physics}, 188(2):1197--1209, August
  2016.

\bibitem{2017JCAP...11..039F}
Luc{\'\i}a {Fonseca de la Bella}, Donough {Regan}, David {Seery}, and Shaun
  {Hotchkiss}.
\newblock {The matter power spectrum in redshift space using effective field
  theory}.
\newblock {\em \jcap}, 2017(11):039, November 2017.

\bibitem{2019MNRAS.483.2078Y}
Victoria {Yankelevich} and Cristiano {Porciani}.
\newblock {Cosmological information in the redshift-space bispectrum}.
\newblock {\em \mnras}, 483(2):2078--2099, February 2019.

\bibitem{2019LRR....22....6H}
Jan {Harms}.
\newblock {Terrestrial gravity fluctuations}.
\newblock {\em Living Reviews in Relativity}, 22(1):6, October 2019.

\bibitem{2018arXiv181202728S}
Zachary {Slepian}.
\newblock {On Decoupling the Integrals of Cosmological Perturbation Theory}.
\newblock {\em arXiv e-prints}, page arXiv:1812.02728, December 2018.

\bibitem{2021MNRAS.501.4004P}
Oliver H.~E. {Philcox}.
\newblock {A faster Fourier transform? Computing small-scale power spectra and
  bispectra for cosmological simulations in $\mathcal{O}(N^2)$ time}.
\newblock {\em \mnras}, 501(3):4004--4034, March 2021.

\bibitem{arfken2013mathematical}
G.B. {Arfken}, H.J. {Weber}, and F.E. {Harris}.
\newblock {\em Mathematical Methods for Physicists: A Comprehensive Guide}.
\newblock Elsevier Science, 2013.

\bibitem{2019arXiv191200065S}
Zachary {Slepian}, Yin {Li}, Marcel {Schmittfull}, and Zvonimir {Vlah}.
\newblock {Rotation method for accelerating multiple-spherical Bessel function
  integrals against a numerical source function}.
\newblock {\em arXiv e-prints}, page arXiv:1912.00065, November 2019.

\bibitem{hamilton2000}
A.~J.~S. {Hamilton}.
\newblock {Uncorrelated modes of the non-linear power spectrum}.
\newblock {\em \mnras}, 312(2):257--284, February 2000.

\bibitem{simonovicfftlog}
Marko {Simonovi{\'c}}, Tobias {Baldauf}, Matias {Zaldarriaga}, John~Joseph
  {Carrasco}, and Juna~A. {Kollmeier}.
\newblock {Cosmological perturbation theory using the FFTLog: formalism and
  connection to QFT loop integrals}.
\newblock {\em \jcap}, 2018(4):030, April 2018.

\bibitem{Dominici_2012}
Diego~E. Dominici, Peter M.~W. Gill, and Taweetham Limpanuparb.
\newblock A remarkable identity involving bessel functions.
\newblock {\em Proceedings of the Royal Society A: Mathematical, Physical and
  Engineering Sciences}, 468(2145):2667–2681, Apr 2012.

\bibitem{2012JChPh.136j4102L}
Taweetham {Limpanuparb}, Joshua~W. {Hollett}, and Peter M.~W. {Gill}.
\newblock {Resolutions of the Coulomb operator. VI. Computation of auxiliary
  integrals}.
\newblock {\em \jcp}, 136(10):104102--104102, March 2012.

\bibitem{2016RSPSA.47260421F}
Bujar~Xh. {Fejzullahu}.
\newblock {On the integral representations for the confluent hypergeometric
  function}.
\newblock {\em Proceedings of the Royal Society of London Series A},
  472(2193):20160421, September 2016.

\bibitem{2019PhRvL.123b0201M}
Satya~N. {Majumdar} and Emmanuel {Trizac}.
\newblock {When Random Walkers Help Solving Intriguing Integrals}.
\newblock {\em \prl}, 123(2):020201, July 2019.

\bibitem{2020JCAP...04..011M}
Azadeh {Moradinezhad Dizgah}, Hayden {Lee}, Marcel {Schmittfull}, and Cora
  {Dvorkin}.
\newblock {Capturing non-Gaussianity of the large-scale structure with weighted
  skew-spectra}.
\newblock {\em \jcap}, 2020(4):011, April 2020.

\bibitem{1994tisp.book.....G}
I.~S. {Gradshteyn} and I.~M. {Ryzhik}.
\newblock {\em {Table of integrals, series and products}}.
\newblock Academic Press, 1994.

\bibitem{abel_1823}
N.H. Abel.
\newblock {\em Solution de quelques problèmes à l'aide d'intégrales
  définies}.
\newblock Cambridge University Press, 1823.

\bibitem{plana_1820}
G.A.A. Plana.
\newblock Sur une nouvelle expression analytique des nombres bernoulliens,
  propre à exprimer en termes finis la formule générale pour la sommation
  des suites.
\newblock {\em Mem. Accad. Sci. Torino}, 25:403–418, 1820.

\bibitem{ButzerPL2011TSFo}
P~L Butzer, P~L Butzer, P~J. S.~G Ferreira, P~J. S.~G Ferreira, G~Schmeisser,
  G~Schmeisser, R~L Stens, and R~L Stens.
\newblock The summation formulae of euler–maclaurin, abel–plana, poisson,
  and their interconnections with the approximate sampling formula of signal
  analysis.
\newblock {\em Resultate der Mathematik}, 59(3):359--400, 2011.

\bibitem{olver1997asymptotics}
F.~Olver.
\newblock {\em Asymptotics and Special Functions}.
\newblock CRC Press, 1997.

\bibitem{saharian2007generalized}
A.~A. Saharian.
\newblock The generalized abel-plana formula with applications to bessel
  functions and casimir effect, 2007.

\bibitem{abramowitz+stegun}
Milton {Abramowitz} and Irene~A. {Stegun}.
\newblock {\em Handbook of Mathematical Functions with Formulas, Graphs, and
  Mathematical Tables}.
\newblock Dover, New York City, 1964.

\bibitem{zygmund}
A.~Zygmund.
\newblock {\em Trigonometric series}, volume I and II combined of {\em 2nd}.
\newblock Cambridge University Press, London and New York, 1969.
\newblock Two volumes bound as one, reprinted with corrections and some
  additions. 2nd edition volume I and II bound together. MR:0236587.

\bibitem{nist_dlmf}
{NIST}.
\newblock {\em NIST Digital Library of Mathematical Functions}.
\newblock DLMF.
\newblock F.~W.~J. Olver, A.~B. {Olde Daalhuis}, D.~W. Lozier, B.~I. Schneider,
  R.~F. Boisvert, C.~W. Clark, B.~R. Miller and B.~V. Saunders, eds.

\end{thebibliography}

\end{document}